\documentclass[12pt]{article}%
\usepackage{amsmath}
\usepackage{amsfonts}
\usepackage{amssymb}
\usepackage{graphicx}%
\newtheorem{theorem}{Theorem}
\newtheorem{acknowledgement}[theorem]{Acknowledgement}

\newenvironment{proof}[1][Proof]{\noindent\textbf{#1.} }{\ \rule{0.5em}{0.5em}}
\begin{document}

\date{}
\title{Some observations on a Kapteyn series}
\author{Diego Dominici \thanks{e-mail: dominicd@newpaltz.edu}\\Department of Mathematics\\State University of New York at New Paltz\\75 S. Manheim Blvd. Suite 9\\New Paltz, NY 12561-2443\\USA\\Phone: (845) 257-2607\\Fax: (845) 257-3571 }
\maketitle

\begin{abstract}
We study the Kapteyn series $%
{\displaystyle\sum\limits_{n=1}^{\infty}}
t^{n}\mathrm{J}_{n}\left(  nz\right)  $. We find a series representation in
powers of $z$ and analyze its radius of convergence.

\end{abstract}


MSC-class: 42C10 (Primary) 30B50, 33C10 (Secondary).

\section{Introduction}

Series of the form
\begin{equation}%
{\displaystyle\sum\limits_{n=1}^{\infty}}
\alpha_{n}\mathrm{J}_{n+\nu}\left[  \left(  n+\nu\right)  z\right]  ,\quad
\nu\in\mathbb{C}\label{KS}%
\end{equation}
where $\mathrm{J}_{n}\left(  \cdot\right)  $ is the Bessel function
of the first kind, are called \textit{Kapteyn series. }The series (\ref{KS}) is convergent and
represents an analytic function throughout the domain in which \cite[17.3]%
{MR1349110}%
\begin{equation}
\Omega\left(  z\right)  <\ \underset{n\rightarrow\infty}{\underline{\lim}%
}\left\vert \alpha_{n}\right\vert ^{-\left(  \nu+n\right)  },\label{ROmega}%
\end{equation}
where%
\begin{equation}
\Omega(z)=\left\vert \frac{z\exp\left[  \sqrt{1-z^{2}}\right]  }%
{1+\sqrt{1-z^{2}}}\right\vert .\label{Omega}%
\end{equation}
Kapteyn series first appeared in Lagrange's solution \cite{lagrange}
\begin{equation}
E=M+2%
{\displaystyle\sum\limits_{n=1}^{\infty}}
\frac{1}{n}\mathrm{J}_{n}\left(  \varepsilon n\right)  \sin\left(  Mn\right)
\label{Kepler}%
\end{equation}
of Kepler's equation \cite{MR1268639}%
\[
M=E-\varepsilon\sin(M).
\]
The solution (\ref{Kepler}) was independently discovered by Friedrich Bessel
in \cite{bessel1}, where he introduced the functions which now bear his name.

Kapteyn series were systematically studied by Willem Kapteyn in his article
\cite{MR911757}, where he proved the following expansion theorem.

\begin{theorem}
Let $f(z)$ be a function which is analytic throughout the region%
\[
D_{a}=\left\{  z\in\mathbb{C}\ |\ \Omega(z)\leq a\right\}  ,
\]
with $a\leq1.$ Then,
\[
f(z)=\alpha_{0}+2\sum\limits_{n=1}^{\infty}\alpha_{n}\mathrm{J}_{n}\left(
nz\right)  ,\quad z\in D_{a}%
\]
where%
\[
\alpha_{n}=\frac{1}{2\pi\mathrm{i}}\oint\Theta_{n}\left(  z\right)  f(z)dz
\]
and the path of integration is the curve on which $\Omega(z)=a.$ The function
$\Theta_{n}\left(  z\right)  $ is the Kapteyn polynomial defined by%
\begin{align*}
\Theta_{0}\left(  z\right)   &  =\frac{1}{z}\\
\Theta_{n}\left(  z\right)   &  =\frac{1}{4}\sum\limits_{k=0}^{\left\lfloor
\frac{n}{2}\right\rfloor }\frac{\left(  n-2k\right)  ^{2}\left(  n-k-1\right)
!}{k!}\left(  \frac{nz}{2}\right)  ^{2k-n},\quad n\geq1,
\end{align*}
where $\left\lfloor \cdot\right\rfloor $ denotes the integer-part function.
\end{theorem}

\begin{proof}
See \cite[17.4]{MR1349110}.
\end{proof}

As a corollary \cite{MR1349110}, one finds that if the Taylor series for
$f(z)$ is
\[
f(z)=\sum\limits_{n=0}^{\infty}a_{n}z^{n}%
\]
then,%
\begin{align}
\alpha_{0}  &  =a_{0}\nonumber\\
\alpha_{n}  &  =\frac{1}{4}\sum\limits_{k=0}^{\left\lfloor \frac{n}%
{2}\right\rfloor }\frac{\left(  n-2k\right)  ^{2}\left(  n-k-1\right)
!}{k!\left(  \frac{n}{2}\right)  ^{n-2k+1}}a_{n-2k},\quad n\geq1.
\label{Inversion}%
\end{align}

In \cite{Kapteyn1}, Kapteyn established the fundamental formula
\begin{equation}
\frac{1}{1-z}=1+2%
{\displaystyle\sum\limits_{n=1}^{\infty}}
\mathrm{J}_{n}\left(  nz\right)  ,\quad\Omega(z)<1, \label{t=1}%
\end{equation}
using the integral representation \cite{MR1349110}%
\[
\mathrm{J}_{n}\left(  nz\right)  =\frac{1}{2\pi i}%
{\displaystyle\int\limits_{\mathcal{H}}}
\left\{  \frac{\exp\left[  \frac{z}{2}\left(  s-\frac{1}{s}\right)  \right]
}{s}\right\}  ^{n}\frac{ds}{s},
\]
where the Hankel contour $\mathcal{H}$ encircles the origin once counterclockwise.

The purpose of this paper is to generalize (\ref{t=1}) by means of the
function%
\begin{equation}
F(z,t)=%
{\displaystyle\sum\limits_{n=1}^{\infty}}
t^{n}\mathrm{J}_{n}\left(  nz\right)  ,\quad z\in\mathbb{C},\quad
t\in\mathbb{R}, \label{F}%
\end{equation}
so that (\ref{t=1}) corresponds to the particular case $t=1.$In Section 2 we
write $F(z,t)$ as a series in powers of $z$%
\begin{equation}
F(z,t)=%
{\displaystyle\sum\limits_{n=1}^{\infty}}
A_{n}(t)z^{n} \label{seriesF}%
\end{equation}
and find the coefficients $A_{n}(t).$ Although one could replace $\alpha
_{n}=t^{n}$ in (\ref{Inversion}) and solve the problem%
\[
t^{n}=\frac{1}{4}\sum\limits_{k=0}^{\left\lfloor \frac{n}{2}\right\rfloor
}\frac{\left(  n-2k\right)  ^{2}\left(  n-k-1\right)  !}{k!\left(  \frac{n}%
{2}\right)  ^{n-2k+1}}A_{n-2k}\left(  t\right)  ,
\]
we take a different approach that only uses the differential equation
satisfied by the Bessel functions. In Section 3 we analyze the radius of
convergence of (\ref{seriesF}) for different ranges of $t.$ Finally, in
Section 4 we solve the inversion problem (\ref{Inversion}) for arbitrary
$\alpha_{n}.$

\section{Power series}

Since the Bessel function $J_{n}(z)$ is a solution of \cite{MR1349110}
\[
z^{2}J_{n}^{\prime\prime}+zJ_{n}^{\prime}+\left(  z^{2}-n^{2}\right)
J_{n}=0,
\]
the function $y_{n}(z)=J_{n}(nz)$ satisfies the ODE%
\begin{equation}
z^{2}y_{n}^{\prime\prime}+zy_{n}^{\prime}+n^{2}\left(  z^{2}-1\right)
y_{n}=0. \label{ODE}%
\end{equation}
We also observe that
\begin{equation}
t^{2}\left(  t^{n}\right)  ^{\prime\prime}+t\left(  t^{n}\right)  ^{\prime
}=n^{2}t^{n}. \label{tn}%
\end{equation}
Using (\ref{ODE}) and (\ref{tn}) in (\ref{F}), we have%
\begin{equation}
z^{2}F_{zz}+zF_{z}=\left(  1-z^{2}\right)  \left(  t^{2}F_{tt}+tF_{t}\right)
, \label{PDEF}%
\end{equation}
where the subscripts denote partial derivatives.

Replacing (\ref{seriesF}) in (\ref{PDEF}), we get the equation%
\begin{equation}
t^{2}A_{n}^{\prime\prime}+tA_{n}^{\prime}-n^{2}A_{n}=t^{2}A_{n-2}%
^{\prime\prime}+tA_{n-2}^{\prime},\quad n\geq1, \label{recu}%
\end{equation}
where we define $A_{n}\equiv0$ for $n=-1,0.$

From (\ref{F}) and (\ref{seriesF}), we have
\[
0=F(z,0)=%
{\displaystyle\sum\limits_{n=1}^{\infty}}
A_{n}(0)z^{n}%
\]
and, from (\ref{t=1}), we see that%
\[%
{\displaystyle\sum\limits_{n=1}^{\infty}}
\frac{1}{2}z^{n}=\frac{z}{2\left(  1-z\right)  }=F(z,1)=%
{\displaystyle\sum\limits_{n=1}^{\infty}}
A_{n}(1)z^{n},
\]
which together imply that
\begin{equation}
A_{n}(0)=0,\quad A_{n}(1)=\frac{1}{2},\quad n\geq1. \label{BC}%
\end{equation}

Solving (\ref{recu}) with the boundary conditions (\ref{BC}), we obtain%
\begin{align*}
A_{1}(t)  &  =\frac{t}{2},\quad A_{2}(t)=\frac{t^{2}}{2},\quad A_{3}%
(t)=\frac{9}{16}t^{3}-\frac{1}{16}t\\
A_{4}(t)  &  =\frac{2}{3}t^{4}-\frac{1}{6}t^{2},\quad A_{5}(t)=\frac{625}%
{768}t^{5}-\frac{81}{256}t^{3}+\frac{1}{384}t,\ldots.
\end{align*}
Thus, the function $A_{n}(t)$ is a polynomial in $t$ of degree $n.$ Writing%
\begin{equation}
A_{n}(t)=%
{\displaystyle\sum\limits_{k=0}^{n}}
C_{k}^{n}t^{k} \label{An}%
\end{equation}
and using (\ref{recu}), we get the recurrence%
\begin{align}
\left(  n^{2}-k^{2}\right)  C_{k}^{n}  &  =-k^{2}C_{k}^{n-2}, & \quad\text{for
}k  &  =0,1,\ldots,n-2\label{Cnk}\\
\left(  n^{2}-k^{2}\right)  C_{k}^{n}  &  =0, & \quad\text{for }k  &
=n-1,\ n.\nonumber
\end{align}
To solve (\ref{Cnk}), we recall that the double factorial function satisfies
\[
\left(  n-k\right)  !!=(n-k)\left(  n-k-2\right)  !!,\quad\left(  n+k\right)
!!=(n+k)\left(  n+k-2\right)  !!,
\]
which gives%
\[
\left(  n^{2}-k^{2}\right)  \frac{1}{\left(  n-k\right)  !!\times\left(
n+k\right)  !!}=\frac{1}{\left(  n-k-2\right)  !!\times\left(  n+k-2\right)
!!}.
\]
This suggests the solution%
\begin{equation}
C_{k}^{n}=\frac{\cos\left(  \frac{n-k}{2}\pi\right)  k^{n}}{\left(
n-k\right)  !!\times\left(  n+k\right)  !!},\quad0\leq k\leq n. \label{Cnksol}%
\end{equation}
Therefore,%
\begin{equation}
A_{n}(t)=%
{\displaystyle\sum\limits_{k=0}^{n}}
\frac{\cos\left(  \frac{n-k}{2}\pi\right)  k^{n}}{\left(  n-k\right)
!!\times\left(  n+k\right)  !!}t^{k} \label{An1}%
\end{equation}
or, after rearranging terms,%
\begin{equation}
A_{n}(t)=\frac{\left(  -1\right)  ^{n}}{n!}%
{\displaystyle\sum\limits_{k=0}^{\left\lfloor \frac{n}{2}\right\rfloor }}
\left(  -1\right)  ^{k}\binom{n}{k}\left(  k-\frac{n}{2}\right)  ^{n}t^{n-2k},
\label{An2}%
\end{equation}

It is clear from (\ref{An1}) that $A_{n}(0)=0.$ To show that $A_{n}%
(1)=\frac{1}{2},$ we use the identity \cite[0.154, \#5]{MR1773820}%
\begin{equation}%
{\displaystyle\sum\limits_{k=0}^{n}}
\left(  -1\right)  ^{k}\binom{n}{k}\left(  k+\alpha\right)  ^{n}=\left(
-1\right)  ^{n}n!, \label{ident.}%
\end{equation}
where $\alpha$ is arbitrary. From (\ref{An2}) we have%
\begin{equation}
A_{n}(t)+A_{n}\left(  \frac{1}{t}\right)  =\frac{\left(  -1\right)  ^{n}}{n!}%
{\displaystyle\sum\limits_{k=0}^{n}}
\left(  -1\right)  ^{k}\binom{n}{k}\left(  k-\frac{n}{2}\right)  ^{n}t^{n-2k}.
\label{An3}%
\end{equation}
Setting $t=1$ in (\ref{An3}), we have%
\[
2A_{n}(1)=\frac{\left(  -1\right)  ^{n}}{n!}%
{\displaystyle\sum\limits_{k=0}^{n}}
\left(  -1\right)  ^{k}\binom{n}{k}\left(  k-\frac{n}{2}\right)  ^{n}%
\]
and from (\ref{ident.}), we conclude that $2A_{n}(1)=1.$

\section{Radius of convergence}

\begin{figure}[ptb]
\begin{center}
\includegraphics[width=4in,angle=270]{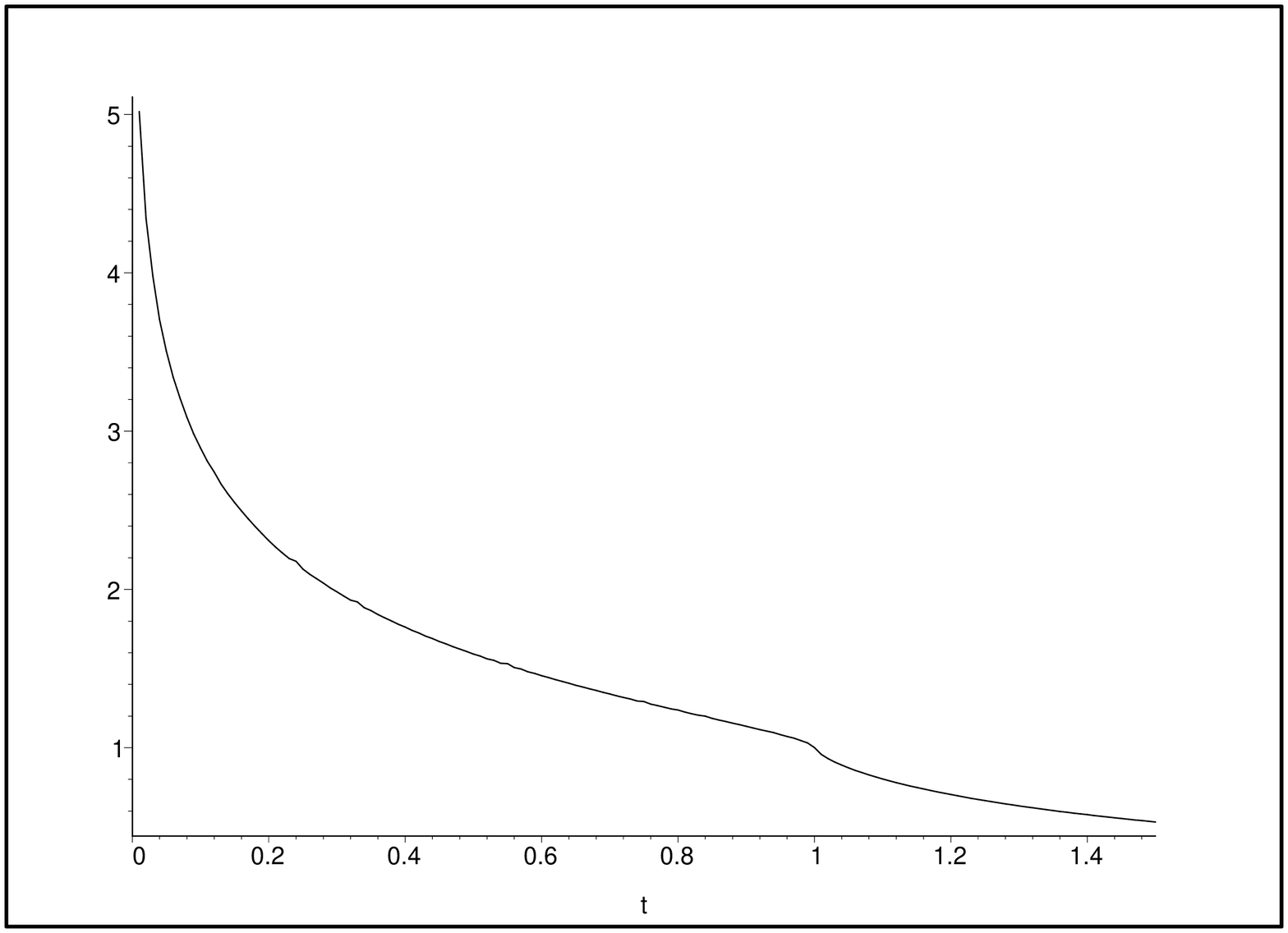}
\end{center}
\caption{A sketch of $\left\vert A_{500}(t)\right\vert ^{-\frac{1}{500}}\simeq
R(t)$.}%
\label{R}%
\end{figure}

We shall now find the radius of convergence $R(t)$ of the power series
(\ref{seriesF}). From (\ref{F}) we observe that $R(t)\rightarrow0$ as
$t\rightarrow\infty$ and $R(t)\rightarrow\infty$ as $t\rightarrow0,$ while
(\ref{t=1}) gives $R(1)=1$ (see Figure \ref{R}). Since from (\ref{An2}) we
have $A_{n}(-t)=\left(  -1\right)  ^{n}A_{n}(t),$ we shall limit our analysis
to $t>0.$

When $t\gg1,$ the largest term in the sum (\ref{An2}) corresponds to $k=0,$
and therefore%

\begin{equation}
A_{n}(t)\sim\frac{1}{n!}\left(  \frac{nt}{2}\right)  ^{n},\quad t\rightarrow
\infty.\label{Antlarge}%
\end{equation}
Using Stirling's formula \cite{AS}%
\[
n!\sim\sqrt{2\pi n}n^{n}e^{-n},\quad n\rightarrow\infty
\]
we obtain%
\begin{equation}
A_{n}(t)\sim\frac{1}{\sqrt{2\pi n}}\left(  \frac{et}{2}\right)  ^{n},\quad
t\rightarrow\infty.\label{Alarge}%
\end{equation}
Thus,%
\begin{equation}
R(t)=\underset{n\rightarrow\infty}{\lim}\left\vert A_{n}(t)\right\vert
^{-\frac{1}{n}}\sim\frac{2e^{-1}}{t},\quad t\rightarrow\infty.\label{Rtlarge}%
\end{equation}
When $t<1,$ (\ref{An2}) is highly oscillatory and there is no approximation
like (\ref{Antlarge}) valid in this range.

\begin{figure}[ptb]
\begin{center}
\includegraphics[width=4in,angle=270]{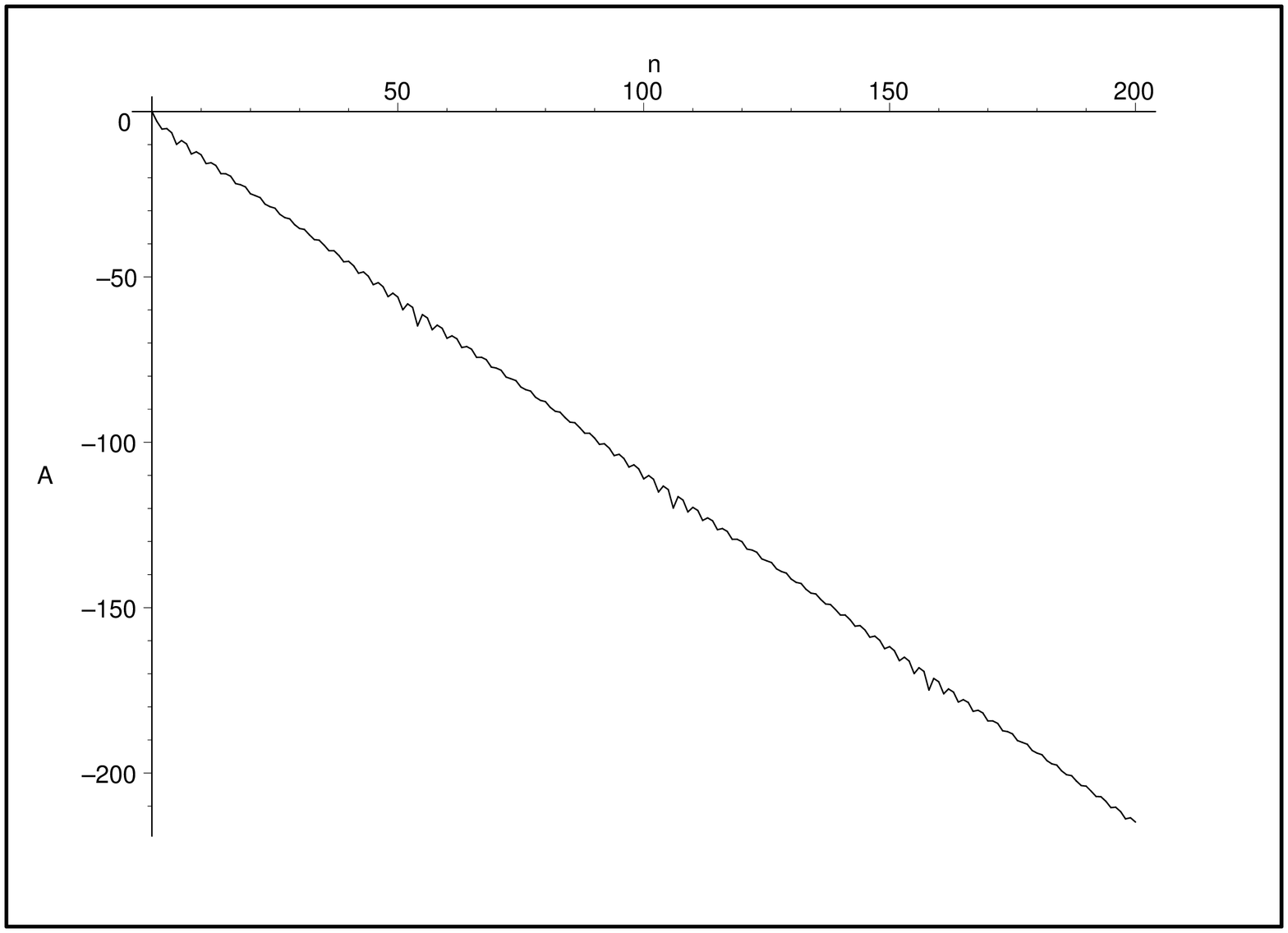}
\end{center}
\caption{A sketch of $\ln\left[\left\vert A_{n}(0.1)\right\vert\right]$.}%
\label{p}%
\end{figure}

To find a formula for $R(t)$ valid for all $t,$ we observe that, for a fixed
value $t>0,$ we have%
\begin{equation}
\ln\left[  \left\vert A_{n}(t)\right\vert \right]  \simeq\beta(t)n,\label{psi}%
\end{equation}
for some function $\beta(t)$ (see Figure \ref{p}). 

Thus, we consider an asymptotic expansion
of the form%
\begin{equation}
A_{n}(t)\sim\frac{1}{2}\theta_{n}(t)\left[  \psi(t)\right]  ^{n},\quad
n\rightarrow\infty,\label{asymA}%
\end{equation}
where $\theta_{n}(t)$ denotes the sign of $A_{n}(t)$ and $\psi(1)=1$. Replacing
(\ref{asymA}) in (\ref{recu}) and setting%
\[
\theta_{n}^{\prime}(t)=\theta_{n}^{\prime\prime}(t)=0,\quad\theta
_{n-2}^{\prime}(t)=\theta_{n-2}^{\prime\prime}(t)=0,
\]
we obtain, to leading order in $n,$%
\begin{equation}
\left(  t\psi^{\prime}\right)  ^{2}\left(  \theta_{n}\psi^{2}-\theta
_{n-2}\right)  -\theta_{n}\psi^{4}=0.\label{eqpsi}%
\end{equation}
Since%
\[
\theta_{n}=\pm1,\quad\theta_{n-2}=\pm1,
\]
and $\psi(t)$ is increasing for $t>0$, we get
\begin{equation}
\psi^{\prime}=\frac{\psi^{2}}{t\sqrt{\psi^{2}\pm1}}.\label{eqpsi2}%
\end{equation}
Solving (\ref{eqpsi2}) subject to $\psi(1)=1,$ we obtain the implicit
solutions%
\begin{equation}
\ln(t)+\frac{\left[  \psi^{2}(t)+1\right]  ^{\frac{3}{2}}}{\psi(t)}%
-\psi(t)\sqrt{\psi^{2}(t)+1}-\operatorname{arcsinh}\left[  \psi(t)\right]
-\sqrt{2}+\ln\left(  \sqrt{2}+1\right)  =0,\label{sol1}%
\end{equation}
for $0<t\leq1$ and%
\begin{equation}
\ln(t)-\frac{\left[  \psi^{2}(t)-1\right]  ^{\frac{3}{2}}}{\psi(t)}%
+\psi(t)\sqrt{\psi^{2}(t)-1}-\ln\left[  \psi(t)+\sqrt{\psi^{2}(t)-1}\right]
=0\label{sol2}%
\end{equation}
for $t\geq1.$

Although we cannot solve (\ref{sol1}) and (\ref{sol2}) exactly, we can
consider the limiting cases as $t\rightarrow0$ and $t\rightarrow\infty.$ Since
$\psi(t)\rightarrow0$ as $t\rightarrow0,$ we obtain from (\ref{sol1})%
\begin{equation}
\psi(t)\sim\frac{1}{-\ln(t)+\sqrt{2}+\ln\left(  \sqrt{2}-1\right)  },\quad
t\rightarrow0.\label{Ap1}%
\end{equation}
When $t\rightarrow\infty,$ we have $\psi(t)\rightarrow\infty,$ and we get from
(\ref{sol2})%
\begin{equation}
\psi(t)\sim\frac{e}{2}t,\quad t\rightarrow\infty,\label{Ap2}%
\end{equation}
which agrees with (\ref{Alarge}). 

\begin{figure}[ptb]
\begin{center}
\includegraphics[width=4in,angle=270]{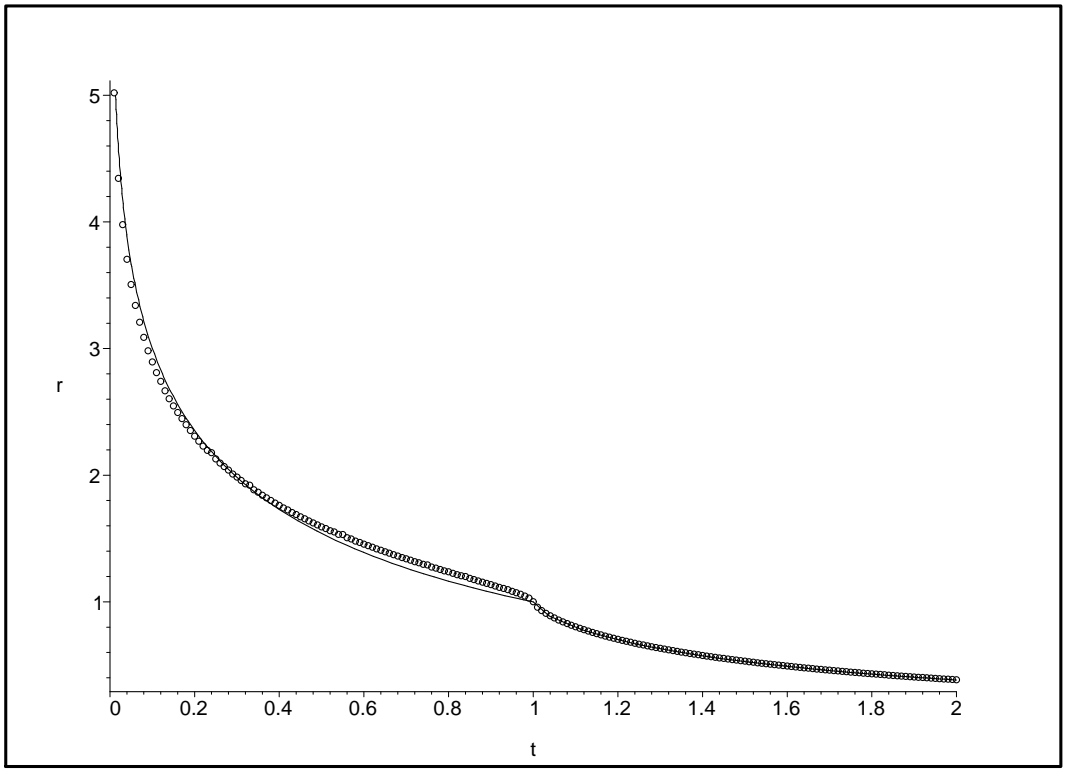}
\end{center}
\caption{A comparison of $R(t)$ (solid curve) and $\left\vert A_{500}(t)\right\vert ^{-\frac{1}{500}}$ (ooo).}%
\label{sol}%
\end{figure}

Since $R(t)=\frac{1}{\psi(t)},$ we have from (\ref{sol1}) and (\ref{sol2}),
after exponentiating%
\begin{equation}
e^{-\sqrt{2}}\left(  1+\sqrt{2}\right)  \frac{R(t)\exp\left[  \sqrt
{1+R^{2}(t)}\right]  }{1+\sqrt{1+R^{2}(t)}}t=1,\quad0<t\leq1,\label{R1}%
\end{equation}
and%
\begin{equation}
\frac{R(t)\exp\left[  \sqrt{1-R^{2}(t)}\right]  }{1+\sqrt{1-R^{2}(t)}%
}t=1,\quad t\geq1.\label{R2}%
\end{equation}
In Figure \ref{sol} we graph the solutions of (\ref{R1}), (\ref{R2}) and the
approximate value of $R(t)$ given by $\left\vert A_{500}(t)\right\vert ^{-\frac{1}{500}}$.

\begin{figure}[ptb]
\begin{center}
\includegraphics[width=4in,angle=270]{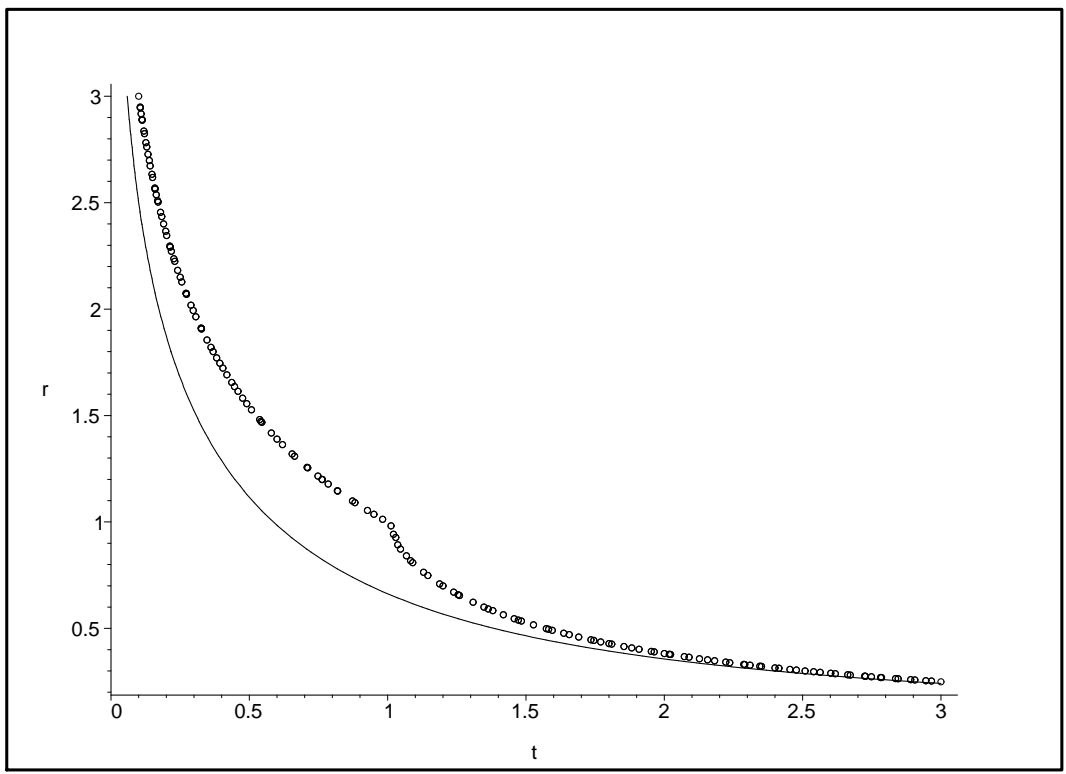}
\end{center}
\caption{A comparison of $r(t)$ (solid curve) and $R(t)$ (ooo).}%
\label{r}%
\end{figure}

Using (\ref{ROmega}) with $\alpha_{n}=t^{n},$ we conclude that the Kapteyn
series (\ref{F}) converges for those $t\in\mathbb{R}$ and $z\in\mathbb{C}$
such that%
\begin{equation}
\left\vert \frac{z\exp\left[  \sqrt{1-z^{2}}\right]  }{1+\sqrt{1-z^{2}}%
}t\right\vert <1.\label{dom1}%
\end{equation}
Replacing $z=re^{\mathrm{i}\omega}$ in (\ref{dom1}), we find that the minimum
value of $r$ corresponds to $\omega=\pm\frac{\pi}{2}.$ Thus, for $t>0,$ the
Kapteyn series (\ref{F}) will converge inside the circle $\left\vert
z\right\vert <r(t),$ with%
\begin{equation}
\frac{r(t)\exp\left[  \sqrt{1+r(t)^{2}}\right]  }{1+\sqrt{1+r(t)^{2}}%
}t=1.\label{dom2}%
\end{equation}
As it was observed in \cite{MR1349110}, the circles of convergence of the
Kapteyn series (\ref{F}) and the power series (\ref{seriesF}) are not equal,
the former being slightly smaller than the latter. However, the radius $r(t)$
given by (\ref{dom2}) and $R(t)$ are asymptotically equal as $t\rightarrow0$
and $t\rightarrow\infty$ (see Figure \ref{r}).

We summarize our results in the following theorem.

\begin{theorem}
The Kapteyn series%
\[
F(z,t)=%
{\displaystyle\sum\limits_{n=1}^{\infty}}
t^{n}\mathrm{J}_{n}\left(  nz\right)  ,
\]
converges for those $t\in\mathbb{R}$ and $z\in\mathbb{C}$ such that%
\[
\left\vert \frac{z\exp\left[  \sqrt{1-z^{2}}\right]  }{1+\sqrt{1-z^{2}}%
}t\right\vert <1.
\]
The function $F(z,t)$ admits the power series representation%
\[
F(z,t)=%
{\displaystyle\sum\limits_{n=1}^{\infty}}
A_{n}(t)z^{n},
\]
with coefficients%
\[
A_{n}(t)=\frac{\left(  -1\right)  ^{n}}{n!}%
{\displaystyle\sum\limits_{k=0}^{\left\lfloor \frac{n}{2}\right\rfloor }}
\left(  -1\right)  ^{k}\binom{n}{k}\left(  k-\frac{n}{2}\right)  ^{n}t^{n-2k}%
\]
and radius of convergence $R(t)$ defined by the implicit equations%
\[
e^{-\sqrt{2}}\left(  1+\sqrt{2}\right)  \frac{R(t)\exp\left[  \sqrt
{1+R^{2}(t)}\right]  }{1+\sqrt{1+R^{2}(t)}}t=1,\quad0<\left\vert t\right\vert
\leq1,
\]%
\[
\frac{R(t)\exp\left[  \sqrt{1-R^{2}(t)}\right]  }{1+\sqrt{1-R^{2}(t)}%
}t=1,\quad\left\vert t\right\vert \geq1.
\]

\end{theorem}

\section{The general problem}

We shall now consider the problem of finding $a_{n}$ in terms of $\alpha_{n},$
where%
\begin{equation}%
{\displaystyle\sum\limits_{n=1}^{\infty}}
a_{n}z^{n}=\sum\limits_{n=1}^{\infty}\alpha_{n}\mathrm{J}_{n}\left(
nz\right)  ,\label{GenK}%
\end{equation}
for arbitrary $a_{n}$ and $\alpha_{n}.$ Rewriting the power series of
$\mathrm{J}_{n}\left(  nz\right)  $ \cite{MR1349110}%
\[
\mathrm{J}_{n}\left(  nz\right)  =\sum\limits_{j=0}^{\infty}\frac{\left(
-1\right)  ^{j}}{j!\left(  n+j\right)  !}\left(  \frac{nz}{2}\right)  ^{n+2j}%
\]
in the form%
\[
\mathrm{J}_{n}\left(  nz\right)  =\left(  \frac{n}{2}\right)  ^{n}%
\sum\limits_{j=0}^{\infty}\frac{\cos\left(  \frac{\pi}{2}j\right)  \left(
\frac{n}{2}\right)  ^{j}}{\left(  \frac{j}{2}\right)  !\left(  n+\frac{j}%
{2}\right)  !}z^{n+j},
\]
we have,%
\begin{align*}
\sum\limits_{n=1}^{\infty}\alpha_{n}\mathrm{J}_{n}\left(  nz\right)    &
=\sum\limits_{n=1}^{\infty}\alpha_{n}\left(  \frac{n}{2}\right)  ^{n}%
\sum\limits_{j=0}^{\infty}\frac{\cos\left(  \frac{\pi}{2}j\right)  \left(
\frac{n}{2}\right)  ^{j}}{\left(  \frac{j}{2}\right)  !\left(  n+\frac{j}%
{2}\right)  !}z^{n+j}\\
& =\sum\limits_{k=1}^{\infty}\sum\limits_{n=1}^{k}\alpha_{n}\left(  \frac
{n}{2}\right)  ^{n}\frac{\cos\left[  \frac{\pi}{2}\left(  k-n\right)  \right]
\left(  \frac{n}{2}\right)  ^{k-n}}{\left(  \frac{k-n}{2}\right)  !\left(
\frac{k+n}{2}\right)  !}z^{k}.
\end{align*}
Therefore, we obtain from (\ref{GenK}) that%
\[
a_{n}=\sum\limits_{n=1}^{k}\alpha_{n}\frac{\cos\left[  \frac{\pi}{2}\left(
k-n\right)  \right]  \left(  \frac{n}{2}\right)  ^{k}}{\left(  \frac{k-n}%
{2}\right)  !\left(  \frac{k+n}{2}\right)  !}%
\]
or%
\begin{equation}
a_{n}=\sum\limits_{n=1}^{k}\alpha_{n}\frac{\cos\left[  \frac{\pi}{2}\left(
k-n\right)  \right]  }{\left(  k-n\right)  !!\left(  k+n\right)  !!}%
n^{k}.\label{GenSol}%
\end{equation}
In particular, setting $\alpha_{n}=t^{n}$ in (\ref{GenSol}), we recover
(\ref{An1}) .

\begin{acknowledgement}
We would like to thank our colleague Michael Adams and the anonymous referees,
for extremely useful comments on earlier versions of this paper.
\end{acknowledgement}

\bibliographystyle{abbrv}
\bibliography{bessel}

\begin{thebibliography}{9}                                                                                                %


\bibitem {AS}M. Abramowitz and I. A. Stegun. \newblock\emph{Handbook of
Mathematical Functions}. \newblock Dover, New York, 9th ed., 1972.

\bibitem {bessel1}F.~W. Bessel. \newblock Analytische {A}ufl\"{o}sung der
{K}eplerschen {A}ufgabe. \newblock {\em Berliner Abh.}, pages 49--55, 1819.

\bibitem {MR1268639}P.~Colwell.
\newblock {\em Solving {K}epler's equation over three centuries}. \newblock
Willmann-Bell Inc., Richmond, VA, 1993.

\bibitem {MR1773820}I.~S. Gradshteyn and I.~M. Ryzhik.
\newblock {\em Table of integrals, series, and products}. \newblock Academic
Press Inc., San Diego, CA, sixth edition, 2000.

\bibitem {MR911757}W. Kapteyn. \newblock Recherches sur les fonctions de
{F}ourier--{B}essel. \newblock {\em Ann. Sci. \'Ecole Norm. Sup.},
x(3):91--120, 1893.

\bibitem {Kapteyn1}W. Kapteyn. \newblock Over Bessel'sche Functi\"{e}n.
\newblock {\em Nieuw Archief voor Wiskunde}, xx:116--127, 1893.

\bibitem {lagrange}J.~L. Lagrange. \newblock Sur le probl\`{e}me de {K}epler.
\newblock {\em Hist. de l'Acad. R. des Sci. de Berlin}, xxv:204--233, 1771.

\bibitem {MR1349110}G.~N. Watson.
\newblock {\em A treatise on the theory of {B}essel functions}. \newblock
Cambridge Mathematical Library. Cambridge University Press, Cambridge, 1995.
\end{thebibliography}

\end{document}